\documentclass{article}

\usepackage{graphicx}
\usepackage[utf8]{inputenc} % allow utf-8 input
\usepackage[T1]{fontenc}    % use 8-bit T1 fonts
\usepackage{hyperref}       % hyperlinks
\usepackage{graphicx}
\usepackage{subcaption}
\usepackage{amsmath,amsfonts,amssymb}
\usepackage{xcolor}
\usepackage[ruled,vlined]{algorithm2e}
\usepackage{floatrow}
\usepackage{wrapfig}
\usepackage{booktabs}
\usepackage{multirow}
\usepackage{array}

\usepackage{tikz}
\usetikzlibrary{shadows,arrows.meta,positioning,backgrounds,fit,chains,scopes}

\tikzset{%
  materia/.style={draw, fill=lightgray!10, text width=20em, text centered, minimum height=1.5em,},
  etape/.style={materia, text width=20em, minimum width=23em, minimum height=7em, rounded corners},
  linepart/.style={draw, thick, color=black!50, -LaTeX, dashed},
  line/.style={draw, thick, color=black!50, -LaTeX},
  ur/.style={draw, text centered, minimum height=0.01em},
  back group/.style={fill=yellow!10,,rounded corners, draw=black!50, dashed, inner xsep=45pt, inner ysep=10pt},
}

\newcommand{\mads}{{MADS}\xspace}
\newcommand{\nomad}{{\sf NOMAD}\xspace}
\newcommand{\hypernomad}{{\sf HyperNOMAD}\xspace}
\newcommand{\hyperopt}{{\sf Hyperopt}\xspace}
\newcommand{\hyperband}{{\sf Hyperband}\xspace}
\newcommand{\pytorch}{{\sf PyTorch}\xspace}

\definecolor{Red}{rgb}{1,0,0}
\definecolor{Blue}{rgb}{0,0.7,0.9}

\newcommand{\bl}{\color{blue}}

\begin{document}
%---------------------------------------------------------------%

%---------------------------------------------------------------%
\title{Use of static surrogates in hyperparameter optimization
\thanks{
\href{https://www.gerad.ca}{\bl GERAD}
  and
\href{https://www.polymtl.ca}{\bl Polytechnique Montr\'eal}
}}
%---------------------------------------------------------------%

%---------------------------------------------------------------%
\author{
    {Dounia Lakhmiri}
  \thanks{\href{mailto:dounia.lakhmiri@polymtl.ca}{\tt \bl dounia.lakhmiri@polymtl.ca}}
  \and
  {S\'ebastien Le~Digabel}
  \thanks{
   \href{https://www.gerad.ca/Sebastien.Le.Digabel/}{\tt \bl www.gerad.ca/Sebastien.Le.Digabel/}
   }
}
%---------------------------------------------------------------%

%\authorrunning{Short form of author list} % if too long for running head

% \institute{Dounia Lakhmiri \at
%             GERAD - Polytechnique Montr\'eal\\
%             \email{dounia.lakhmiri@polymtl.ca}
%         \and S\'ebastien Le~Digabel \at
%             GERAD - Polytechnique Montr\'eal\\
%             \email{sebastien.le.digabel@gerad.ca}}

% \date{Received: date / Accepted: date}
% The correct dates will be entered by the editor

\maketitle

%---------------------------------------------------------------%
\begin{abstract}
Optimizing the hyperparameters and architecture of a neural network is a long yet necessary phase in the development of any new application. This consuming process can benefit from the elaboration of strategies designed to quickly discard low quality configurations and focus on more promising candidates. This work aims  at enhancing \hypernomad, a library that adapts a direct search derivative-free optimization algorithm to tune both the architecture and the training of a neural network simultaneously, by targeting two keys steps of its execution and exploiting cheap approximations in the form of static surrogates to trigger the early stopping of the evaluation of a configuration and the ranking of pools of candidates. These additions to \hypernomad are shown to improve on its resources consumption without harming the quality of the proposed solutions.\\
\end{abstract}
%---------------------------------------------------------------%

%---------------------------------------------------------------%
\section{Introduction}
%---------------------------------------------------------------%

The efficacy of deep neural networks (DNN) to discover patterns in complex datasets is seen throughout multiple applications where the appropriate variant of a DNN often manages to score higher than human experts or other machine learning algorithms and even to push the current state of the art. A particular class of DNN includes convolutional neural networks (CNN) which are used for image classification, image segmentation, or object detection. They have been gaining in popularity and attention from the scientific community in the recent years as they are at heart of important technological advances such as imitation learning~\cite{10.1007/978-3-030-61377-8_26} or medical imagery analysis~\cite{marques2020automated} to name a few. 
    
An important challenge when adopting a deep learning approach for a new task is to find the appropriate neural network by deciding on a set of hyperparameters that determine the architecture, i.e, the number of layers, type of blocks and connection, and the training regime. The final score of the network is highly sensitive to the tuning of these hyperparameters and there is no rigorous or intuitive rule that directly provides a range for their optimum values. This hyperparameter optimization (HPO) problem can therefore be framed as a blackbox optimization problem where the objective value is the result of a computation with no analytical formula and no known derivatives. Moreover, this process is time consuming and expensive in computational resources as each configuration trial is equivalent to training a new network long enough to infer its generalization score. In a supervised learning setting, the HPO problem can be expressed as
\begin{equation}
    \label{eq:problem}
    \min_\Phi  f(\Phi, \theta^*) \quad 
 \text{ with }   \theta^* \in \operatorname*{argmin}_\theta \mathcal{L}_{\Phi, \theta}(X,Y)
\end{equation}
where $f$ is the objective function that corresponds to a measure of performance of the neural network, usually the validation loss, $\Phi$ is the mixed integer vector of all the hyperparameters that define the network, $X,Y$ are the validation data and labels, $\mathcal{L}$ is the loss function that the network optimizes during the training by updating the weights $\theta$.

Derivative-free optimization (DFO)~\cite{AuHa2017,CoScVibook} provides a class of methods that are well suited to tackle such blackbox HPO problems as they do not need the explicit expression of the objective function and/or the constraints, nor do they rely on the derivatives for their execution. Two classes of DFO algorithms can be defined: Model-based and direct search methods. Model-based algorithms use a static or dynamic surrogate function $\hat{f}$ as an approximation of the true objective $f$ to guide the optimization. Static surrogates, also called simplified physics surrogates, are defined before the start of the optimization and remain unchanged during the execution in contrast to dynamic surrogates, such as Gaussian processes that are updated with each iteration of the DFO method or with each new evaluation of the objective function. Direct search methods base their exploration of the search space solely on the values of the evaluated points and explore said space through a pattern as it is the case for the Mesh Adaptive Direct Search (\mads)~\cite{AuDe2006} or on a simplex in the case of the Nelder-Mead algorithm~\cite{NeMe65a}. The properties of DFO methods explain their popularity in the context of HPO of deep neural networks where they are often included in specialized libraries such as \hyperopt~\cite{bergstra2013making} or {\sf Or\'ion}~\cite{xavier_bouthillier_2019_3478593}. Similarly, the \hypernomad toolbox~\cite{hypernomad,hypernomad_paper} is developed as an adaptation of \mads to simultaneously optimize the architecture and the training phase of a CNN for a given dataset as expressed
in~\eqref{eq:problem}. This open source library was shown to have a competitive performance against other popular approaches such as Bayesian optimization~\cite{bergstra2011algorithms} or a random search~\cite{bergstra2012random} on the MNIST~\cite{mnist}, Fashion-MNIST~\cite{xiao2017fashion} and CIFAR-10~\cite{krizhevsky2009learning} datasets.

The objective of this work is to develop a set of protocols that speed up the HPO process in \hypernomad. The proposed measures are based on the use of two static surrogates that affect different steps of the algorithm execution: first the training log of the best encountered network serves as a baseline of comparison with the currently evaluated point and allows for early stopping decisions if the performance of the current network seems unsatisfactory. Second, another static surrogate is employed to rank the candidates to be evaluated at each iteration of \hypernomad which applies an opportunistic evaluation strategy, meaning that an iteration is interrupted as soon as a better score is recorded, therefore disregarding the other pool of candidates. The implementation of these two surrogates is shown to significantly speed up the HPO process thus saving expensive resources. 

The rest of the document is structured as follows: Section~\ref{sec:literature} provides an overview of the literature regarding strategies for early stopping and accelerating the costly HPO problem. Section~\ref{sec:hypernomad} goes through the algorithm behind \hypernomad to highlight where the added strategies are implemented. Section~\ref{sec:surrogates} dives into the details of the proposed speed ups and the management of each surrogate. Finally, Section~\ref{sec:results} compares the original version of \hypernomad with the proposed additions and a synthesis of these results is presented in the discussion. 

%---------------------------------------------------------------%
\section{Related work}\label{sec:literature}
%---------------------------------------------------------------%

% \begin{itemize}
%     \item Early stopping, but when ? presents a number of criterias to analyse the training of a single machine learning algorithm that measure the improvement of the validation error, the stagnation of that error... 
%     \item Hyperband has a similar management to early stopping.
%     \item NAS and AutoML surveys, the estimation of score part. 
% \end{itemize}

The common problem of fine tuning a neural network to yield the best performance measure for a specific dataset is highly complex and consuming in terms of time and computational resources, so much so that it is often treated in two separate steps: one for finding the neural architecture and one for optimizing the hyperparameters related to the training of the network. Regardless of which phase is considered, the problem can be formulated as an expensive blackbox optimization problem~\cite{AuHa2017,CoScVibook} since the evaluation of the objective function is equivalent to training a new neural network with a specific configuration for a certain amount of time in order to estimate its final accuracy. That necessary amount of time is unknown a priori and has a direct impact on the quantity of resources needed to carry out this task. Moreover, experience shows that not every hyperparameter configuration can yield satisfactory results and it is best to develop the tools that quickly detect such cases and prematurely stop their training, thus saving valuable resources than can be allocated to more promising candidates. 

Early stopping is first introduced as a regularization technique to ovoid the overfitting issue that arises in machine learning. Theoretically, the training of a network, and more generally any machine learning algorithm, should be interrupted as soon as the validation error starts to increase as that behaviour indicates a training saturation and an overfitting to the training data. However, as the authors of~\cite{Prechelt2012} point out, this criteria can not be directly applied on a real validation error curve that tend to be irregular
% as shown in Figure~\ref{fig:real_val_curve}. 
This work also provides a set of adapted criteria that prompt early stopping when the increase of the validation error exceeds a predetermined threshold, when the deterioration of the validation error is greater than the improvement of the training error or when the generalisation error consistently increases over a certain number of {\em epochs}. An epoch describes one pass on the entire training data. The training usually consists in multiple epochs before reaching convergence. Each of these criteria is shown to provide good a trade-off between the overall training time and the final performance of the network on a collection of problems. 
Early stopping is also exploited in \hyperband~\cite{li2017hyperband} which adapts the random search (RS)~\cite{bergstra2012random} with an efficient resource management scheme. This algorithm starts with allocating a fraction of the resources to each new configuration and compares the resulting validation loss, half of the the low performing candidates are stopped and the other half is granted more resources and the process is repeated throughout the execution of \hyperband with only the relative top performers being trained with the full resources.
%as illustrated in Figure~\ref{fig:hyperband}. 
Compared to other HPO algorithms such as Bayesian methods, \hyperband allows for a significant speed up on multiple datasets such as MNIST and CIFAR-10. The same principle of ``starting many, stopping early, continuing some'' is applied in~\cite{dodge2020fine} where the authors obtained the best results when discarding the least performing networks after $20\%-30\%$ of the total training budget. All aforementioned works rely on the observed validation curves, yet they highlight a different aspect of the early stopping decisions. In~\cite{Prechelt2012}, stopping the training is only dependent on the scores of the network itself whereas~\cite{dodge2020fine,li2017hyperband} stop networks for their relative scores compared to a pool of candidates. 

Besides, early stopping strategies can base their decisions on the estimation of future scores instead of the observed ones only. Survey~\cite{elsken2018neural} dedicates a section to some of the methods used to estimate the final validation performance, or to a lesser extent, the expected scores in future iterations. These strategies extrapolate the validation or
training curves~\cite{domhan2015speeding} via a regression method by factoring the observed scores and, possibly, the architectural and learning hyperparameters of the network~\cite{baker2017accelerating,Klein2017LearningCP}. The usefulness of these estimators in speeding-up hyperparameter optimization or neural architecture search is shown in each work but their main challenge is their need to be trained on a substantial amount of data before having reliable predictions. Another score estimation route uses low fidelity surrogates such as training a network with a lower epoch budget~\cite{zela2018automated}, on a subset of the dataset~\cite{klein2017fast,trofimov2020multifidelity} or on lower resolution images~\cite{chrabaszcz2017downsampled}.
The predicted scores with these methods are expected to be under-estimations of the real validation scores, however they can still be used as a ranking tool to separate the candidate configurations between the top and low performers as long as a good balance in the low fidelity estimate is found so that the under-estimations are not too large.
% to be unexploitable.  

% \begin{figure}
%     \begin{subfigure}[b]{0.43\textwidth}
%          \centering
%          \includegraphics[width=\textwidth]{figs/real_val_loss.pdf}
%          \caption{Example of a real validation error curve.}
%          \label{fig:real_val_curve}
%      \end{subfigure}
%      \hfill
%      % \hspace{5cm}
%      \begin{subfigure}[b]{0.43\textwidth}
%          \centering
%          \includegraphics[width=\textwidth]{figs/hyperband.png}
%          \caption{Example of Hyperband's allocation of training resources. Image taken~\cite{li2017hyperband}.}
%          \label{fig:hyperband}
%      \end{subfigure}
%      \caption{\rd New figures à produire}
% \end{figure}

%---------------------------------------------------------------%
\section[The HyperNOMAD package]{The \hypernomad package}
\label{sec:hypernomad}
%---------------------------------------------------------------%

% Quick review of HyperNOMAD and MADS to explain where the surrogates will be used.

\begin{figure}[h]
%   \begin{center}
%     \includegraphics[scale=0.65]{figs/hypernomad_acc.pdf}
%   \end{center}
\begin{tikzpicture}
     [
    start chain=p going below,
    every on chain/.append style={etape},
    every join/.append style={line},
    node distance=1 and -.25,
  ]
  {
    \node [on chain, join] {\textbf{Initial parameters}: dataset, initial point, budget evaluation};
    
    % Grande boite grise optimizer
    \node (optim) [on chain, join] {};
        \node at ([yshift=-1em] optim.north){\textbf{HyperNOMAD optimizer}};
        \node[draw, thick, shape=rectangle, minimum width=3.5cm, minimum height=1.5cm, rounded corners] at ([xshift=6em, yshift=-1em]optim.west){\begin{tabular}{@{}c@{}}Hyperparameters,\\block structure\\ Neighbors\end{tabular}};
        \node (nomad) [draw, thick, shape=rectangle,  minimum width=3.5cm, minimum height=1.5cm, rounded corners] at ([xshift=-6em, yshift=-1em]optim.east){NOMAD} ;

    % Grande boite grise BB
    \node (bb) [on chain, thick] {};
        \node at ([yshift=-1em] bb.north){\textbf{Blackbox}};
        \node[draw, thick, shape=rectangle, minimum width=3.5cm, minimum height=1.5cm, rounded corners] at ([xshift=6em, yshift=-1em]bb.west){Construct network};
        \node[draw, thick, shape=rectangle,  minimum width=3.5cm, minimum height=1.5cm, rounded corners] at ([xshift=-6em, yshift=-1em]bb.east){\begin{tabular}{@{}c@{}}Network training,\\validation, testing\end{tabular}
                            };

    \node (tmp1) [right=1cm of optim] {};
    \draw[-] (optim.east) -- (tmp1.center);
    \draw[->] (tmp1.center) |-  (bb.east) {};
    
    \node[right=-3cm of optim, below= 1.5cm of optim.east]{New point};
    
    \node (tmp2) [left=1cm of bb] {};
    \draw[-] (bb.west) -- (tmp2.center);
    \draw[->] (tmp2.center) |-  (optim.west) {};
    
    \node[xshift=1.5em, above= 1.5cm of bb.west]{Validation accuracy};
  }

    \begin{scope}[on background layer]
        \node (bk1) [back group] [fit=(p-2) (p-3)] {};
        \node [draw, thick, yellow!5!green!10, fill=yellow!5!green!10, rounded corners, fit=(p-1), inner xsep=15pt, inner ysep=10pt] {};
    \end{scope}

  \end{tikzpicture}

  \caption{The \hypernomad workflow is represented by the communication of its two components. The optimizer suggests new candidates based on the \nomad software and the blackbox trains the corresponding neural network to return its validation or test accuracy as the objective function value. Image adapted from~\cite{hypernomad_paper}.}
  \label{fig:HyperNOMAD}
\end{figure}
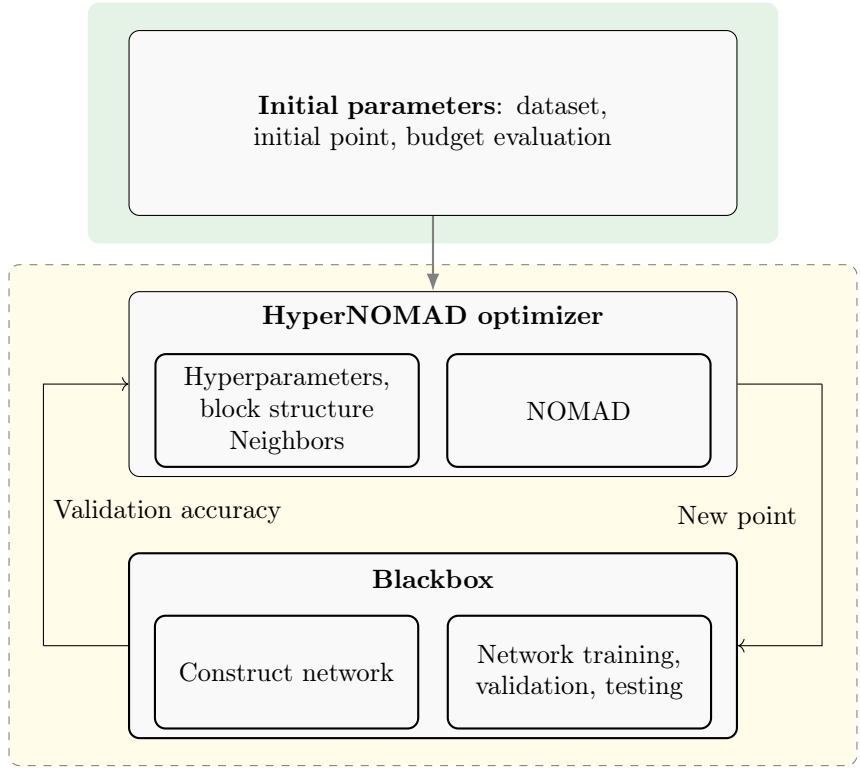

The open-source library \hypernomad~\cite{hypernomad,hypernomad_paper} is designed as an adaptation of the \nomad software~\cite{Le09b} to optimize the hyperparameters of deep neural networks as formulated in~\eqref{eq:problem}.  This package allows searching for both the architecture and the convolutional network's training regime for a specific dataset. It contains two main components: the blackbox and the optimizer, as illustrated in Figure~\ref{fig:HyperNOMAD} and presented in the following section.

%-------------------------%
\subsection{The blackbox}
%-------------------------%
This component regroups the {\sf Python} and {\sf Pytorch}~\cite{paszke2019pytorch} modules responsible for creating the equivalent deep neural network to the provided vector of hyperparameters. The blackbox module fully trains the resulting network before returning its validation accuracy or its validation loss as a performance measure. The package provides some computer vision datasets such as MNIST~\cite{mnist}, Fashion-MNIST~\cite{xiao2017fashion}, CIFAR-10/100~\cite{krizhevsky2009learning} and STL10~\cite{coates2011analysis} but also allows to plug-in a new dataset if needed.
% In the original version of \hypernomad, each network is trained during $200$ epochs, and early stopping can occur if the standard variation of the validation loss for the last $50$ epochs is lower than $10^{-3}$.\\
\subsection{The optimizer} 
The second module is first responsible for handling the collected categorical, integer, or real HPs. Categorical HPs include non-ordinal variables such as the activation function or the choice of the training optimizer and some ordinal variables such as the number of convolutional or fully connected layers. Then, the module launches the optimization through the \nomad software, which is the official implementation of the Mesh Adaptive Direct Search (\mads)~\cite{AuDe2006} algorithm further discussed in Section~\ref{sec:mads}.

In \hypernomad, each convolutional layer is defined by five HPs: the number of output channels, the kernel size, the stride, the padding, and the pooling size; And each fully connected layer is determined by the number of nodes it contains. Hence, if $n_1$ indicates the number of convolutional layers and $n_2$ the number of fully connected layers, the network's architecture is determined by $5n_1 + n_2$ HPs. Considering the remaining HPs such as the optimizer, batch size or dropout rate, the dimension of the complete HPO problem becomes $5n_1+n_2+10$. This dimension varies during the optimization if the values of $n_1$ and $n_2$ evolve, explaining why these HPs are considered categorical variables. A change in their values signals a change of search space, which \hypernomad handles through the extended polls in Algorithm~\ref{algo-mads} with a predefined neighborhood structure that states that each CNN has five neighbors.
Four of these neighbors are found by adding/subtracting one convolutional/fully connected layer and the fifth neighbor keeps the same architecture and changes the optimizer for the training task. 

%---------------------------------------------------------------%
\subsection{The mesh adaptive direct search (\mads) algorithm}
\label{sec:mads}
%---------------------------------------------------------------%

As previously stated, the optimizer in \hypernomad launches the \mads algorithm to explore the HPs search space. At each iteration $k$, \mads works on a mesh $M_k$ defined by a mesh size $\Delta_k$ that is expanded or reduced depending on whether the previous iteration is successful or not. An iteration is successful if an improvement on the incumbent is recorded; otherwise, it is a failure. 

Each iteration is composed of two steps. First, the {\em search} step is an optional and flexible phase containing any exploration strategy as long as a finite number of trials projected on the mesh $M_k$ are evaluated. Then the {\em poll} step starts by listing a set of poll points around the incumbent $P_k = \{x_k + \Delta_k d \mid d \in D_k\}$ where $D_k$ is a set of search directions that form a positive basis. The search directions are generated so that the equivalent unit directions become asymptotically dense in the unit sphere. 
The points $p_k\in P_k$ are evaluated with an {\em opportunistic} strategy, which means the poll step is interrupted as soon as a configuration scores better than the incumbent $x_k$, even if other poll points are still available and hence will not be tested. Algorithm~\ref{algo-mads} summarizes the key steps of the \mads algorithm with an emphasis on where surrogate functions can be integrated. 

In the general framework, surrogate functions are optionally plugged into the \mads execution at different stages. For example, a static or dynamic surrogate can be used to explore a wide range of the space during the search phase to suggest new candidates. 
The poll step can also benefit from surrogate assistance in few distinct aspects, such as deciding on the evaluations to interrupt or providing a ranking so that \mads evaluates the most promising candidate first. Such rankings coupled with the opportunistic strategy are essential to target better solutions faster~\cite{loic}, especially when time and resource constraints are as severe as in the current context. The present work focuses mainly on improving the poll step by targeting both the resource allocation and the ranking aspects of poll points evaluation. The incorporated improvements are discussed in detail in Section~\ref{sec:surrogates}.

% Algorithm~\ref{algo-mads}. Each iteration $k$ defines a mesh $M_k=\{x + \Delta_k^{m}Dz, z \in \mathbb{N}^{n_D}, x \in C\}$ where $C$ regroups the points evaluated thus far, $D$ is a matrix with columns that define a positive spanning set and $\Delta_k$ is the current mesh size. If $x_k$ is the best feasible solution at the current iteration, a pool of candidates $P_k$ is defined during the poll step as $P_k = \{x_k + \Delta_k^{m}d \mid d \in D_k\}$. Originally, the new candidates are evaluated following the direction of the last success first {\rd plus, an opportunistic strategy is applied, meaning that} the poll step is interrupted as soon as a configuration scores better than the incumbent $x_k$. 

%---------------------------------------------------------------%
\begin{algorithm}[ht]
        \caption{\mads with static surrogates for ranking and early stopping.}
        \begin{flushleft}
        \textbf{[0] Initialization}\\
        \hspace*{5mm} iteration $k = 0$, configuration $x_0$, mesh size $\Delta_0$ and mesh $M_0$,\\
        \hspace*{5mm} ranking surrogate $S_1$ and early stopping surrogate $S_2$\\
        
        \textbf{[1] Search of iteration $k$ (optional)}\\
        \hspace*{5mm} Construct a set of mesh points and evaluate them \\
        \hspace*{5mm} If there is a success, go to \textbf{[3]}

        \textbf{[2] Poll of iteration $k$}\\
        \hspace*{5mm} Define the poll set $P_k$ of new candidates around $x_k$,\\
        \hspace*{5mm} along search directions that define a positive basis.\\
        \hspace*{5mm} {\em Sort the poll points with surrogate $S_1$} \\
        \hspace*{5mm} {\em Evaluate the points in $P_k$ with possible early stopping ($S_2$)}\\
        \hspace*{5mm} If there is a success, interrupt the evaluations and go to \textbf{[3]}

        \textbf{[3] Updates of iteration $k$}\\
        \hspace*{5mm} Update $\Delta_k, x_k, M_k$ depending on the success of the previous phases\\
        \hspace*{5mm} If no stopping condition is satisfied: $k \leftarrow k+1$ and go to \textbf{[1]}
        \end{flushleft}
    \label{algo-mads}
\end{algorithm}
%---------------------------------------------------------------%

%---------------------------------------------------------------%
\section{Proposed approach}\label{sec:surrogates}
%---------------------------------------------------------------%

This section details the proposed enhancements to speed-up the resolution of the HPO problem with \hypernomad. This objective is achieved by introducing an early stopping mechanism that quickly interrupts the poor performing candidates, coupled with a ranking strategy that allows to evaluate the most promising poll candidates first. All tests in this section are on the MNIST dataset and start the HPO process from two configurations. First with the default initialization of \hypernomad, noted $p_1$, that corresponds to a network with one convolutional layer and two fully connected layers, which amounts to 17 hyperparameters. The second initialization, noted $p_2$, adds a convolutional layer to the default point, which amounts to 22 hyperparameters.

%---------------------------------------------------------------%
\subsection{Early stopping}
%---------------------------------------------------------------%

Implementing an adequate early stopping strategy is paramount to the efficient execution of any hyperparameter optimization technique as the complete evaluation of each encountered candidate is merely unreasonable due to the number of wasted resources this direct approach can cause. The challenge is to detect when a network is worth training further and when it needs to be interrupted based on the observed validation curve and compared to other candidates previously trained. 

\begin{figure}[ht]
    \begin{subfigure}[t]{0.49\textwidth}
        \centering
        \includegraphics[width=1.1\textwidth]{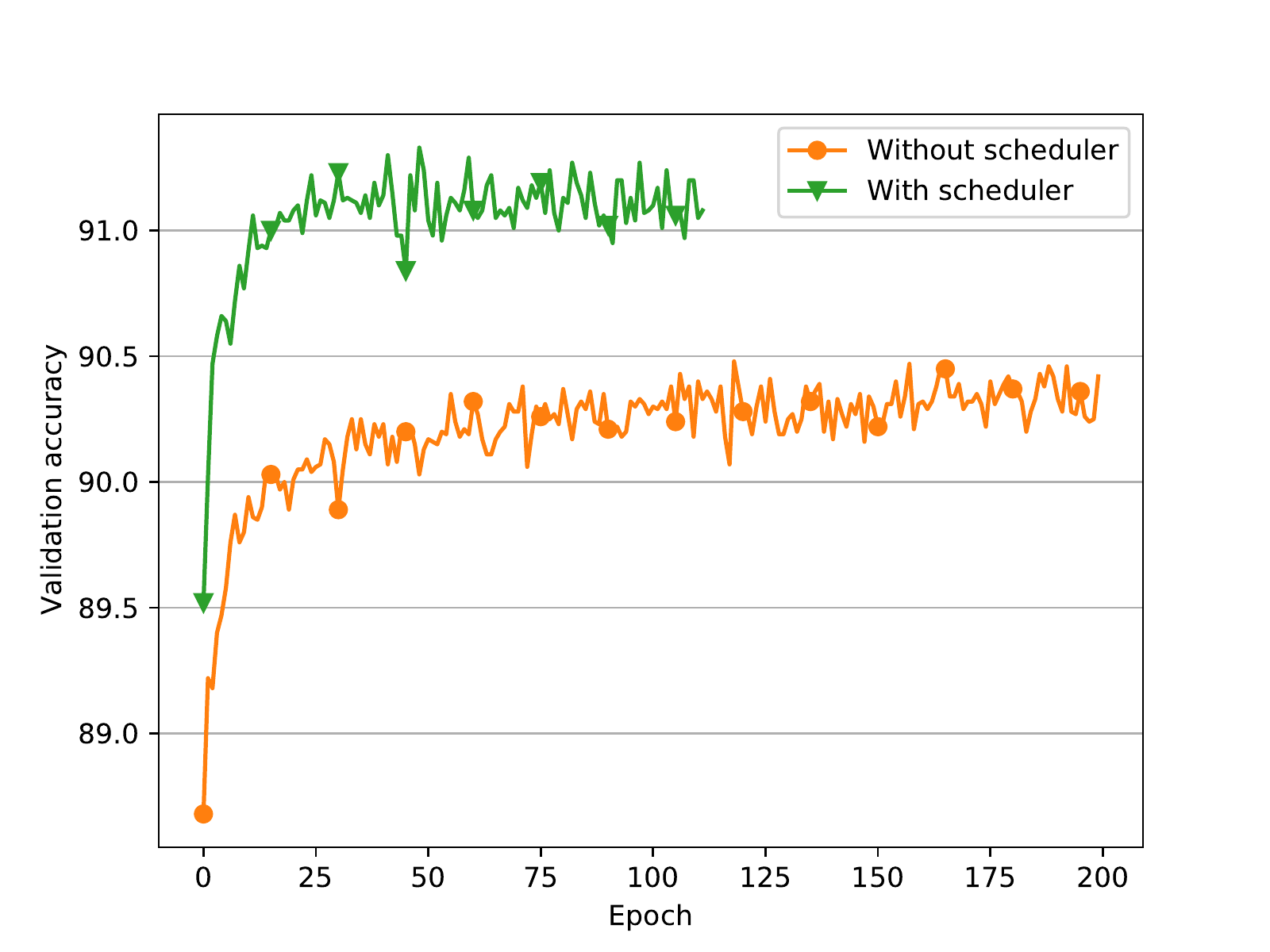}
        \caption{Example of early stopping due to the scheduler.}
        \label{fig:scheduler}
    \end{subfigure}
    \hfill
    \begin{subfigure}[t]{0.49\textwidth}
        \centering
        \includegraphics[width=1.1\textwidth]{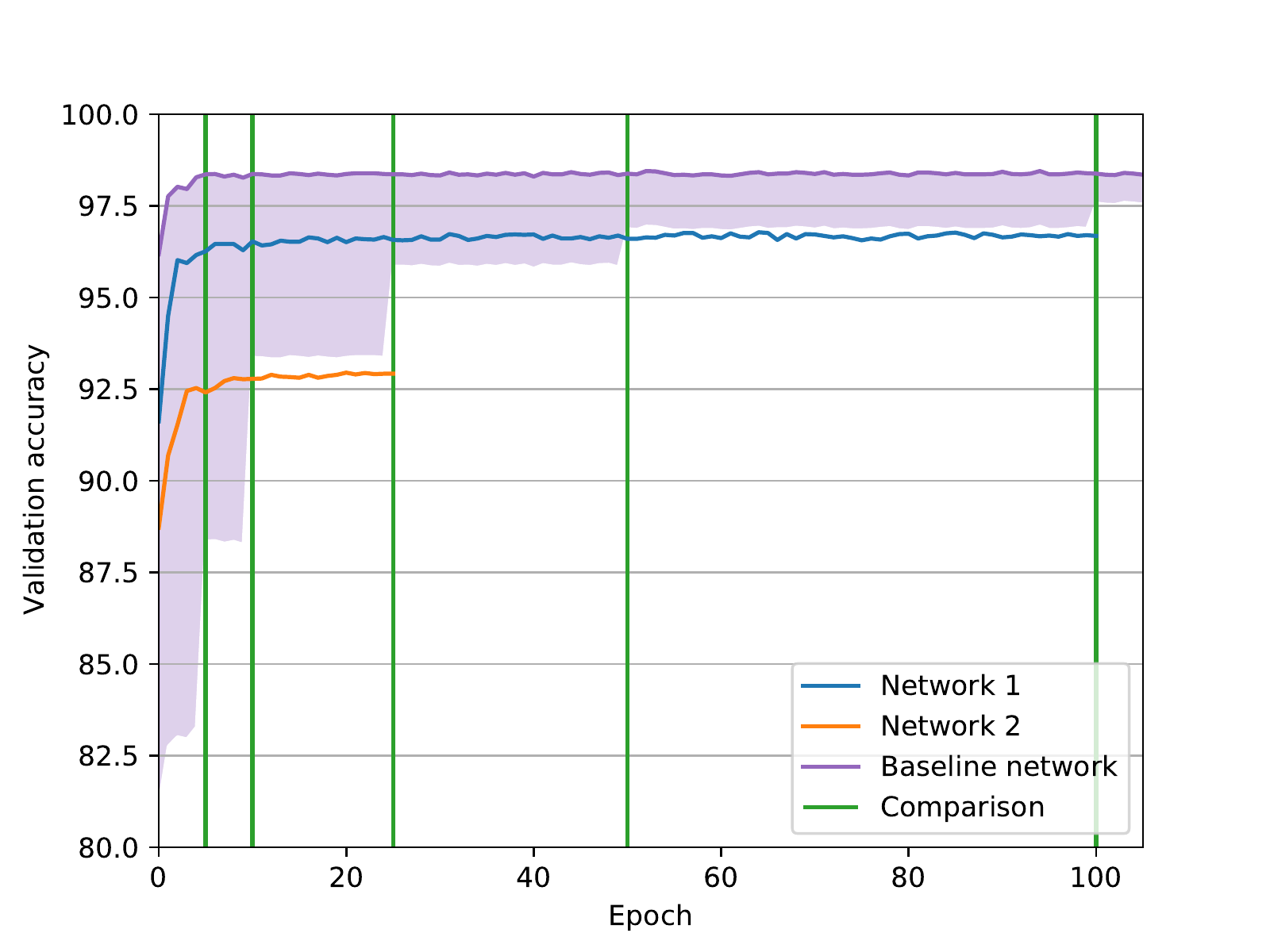}
        \caption{Example of early stopping due to the comparison with the baseline network.}
        \label{fig:baseline}
     \end{subfigure}
    \caption{Two early stopping criteria: a scheduler that detects a plateau and reduces the learning rate until it hits the lower limits (left) and the comparison with a baseline network that stops any evaluation outside the envelope defined by the relative errors compared to the baseline scores(right).} 
    \label{fig:earystopping}
\end{figure}

The validation accuracy of a network should gradually improve during the training of a DNN until hitting a plateau or a maximum value after which the accuracy decreases.  Such scenarios need to be detected to avoid depleting the remaining training budget on a configuration that can not improve its current validation score any further. However, a plateau does not necessarily mean that the best validation accuracy is reached but rather that the current training HPs, especially the learning rate, need to be updated. To this effect, the \pytorch library~\cite{paszke2019pytorch} provides multiple scheduler variants to manage the learning rate during the training. In particular, the \texttt{ReduceLROnPlateau} scheduler can detect plateaus and prompt the desired change in the learning rate. A common strategy is to reduce the learning rate progressively until it reaches a predefined minimum value, which triggers an early stopping. Figure~\ref{fig:scheduler} illustrates the benefit of such a mechanism as adding the scheduler allows saving half the training budget in that particular case. 

Early stopping is also a decision based on the relative score of the current network compared to the candidates previously evaluated. Such comparison targets the networks whose validation scores are too far from the best solution seen thus far, even if that score gradually improves as the training advances. When training a new candidate, its validation accuracy curve is compared against the best scoring network, called the baseline, at specific epochs: $[5, 10, 25, 50, 100, 125, 150 ]$ with each milestone having an associated error margin: $[0.5, 0.6, 0.7, 0.8, 0.85, 0.9, 0.95]$. After 5 epochs, a new configuration needs to score at least 50\% of the baseline score, 60\% at 10 epochs, and so on. As such, the error margins define an envelope under the baseline curve, and any network that scores lower than the allowed error margin is interrupted at the next milestone epoch. This early stopping mechanism is illustrated in Figure~\ref{fig:baseline}.
If a better solution is found, the corresponding network becomes the new baseline, and the envelope is redefined with its validation accuracy curve. Consequently, as the HPO advances, better baselines are found, and only the high scoring networks are allowed a high training budget and the low performances are quickly detected, and interrupted. 

Early stopping effects are first tested on the MNIST dataset by comparing the original version of \hypernomad with a variant that implements one early stopping strategy. Every HPO execution starts from the same initial configuration and allows for 200 blackbox evaluations (BBE), meaning that 200 different configurations are trained and tested. The comparison aims at observing the effect on the score of the best solution obtained and the quantity of resources needed to carry out the entire optimization task. Resource consumption is measured in terms of overall wall-clock time and total number of epochs used in training every visited configuration. In this phase, every execution is run on one Nvidia P100 GPU, and the batch size is fixed to 256.
The tested early stopping strategies are:
\begin{itemize}
    \item Default settings of \hypernomad: stops if the validation accuracy does not surpass 12\% after 25 epochs or when the standard deviation of the validation loss over the last 50 epochs is lower than $10^{-3}$.
    \item Last success:  stops the training when the last improvement on the validation accuracy is recorded more than $25$ epochs ago.
    \item Scheduler alone:  corresponds to the \texttt{ReduceLROnPlateau} scheduler that divides by $10$ the learning rate if the validation accuracy has not improve for the last 25 epochs until the learning rate becomes lower than $10^{-8}$. 
    \item Baseline and scheduler: takes the previous scheduler scheme and adds the relative scores comparison by defining a baseline and the corresponding envelop as previously detailed.
\end{itemize}

\begin{figure}[ht]
    \centering
    \begin{subfigure}[t]{0.49\textwidth}
        \centering
        \includegraphics[width=\textwidth]{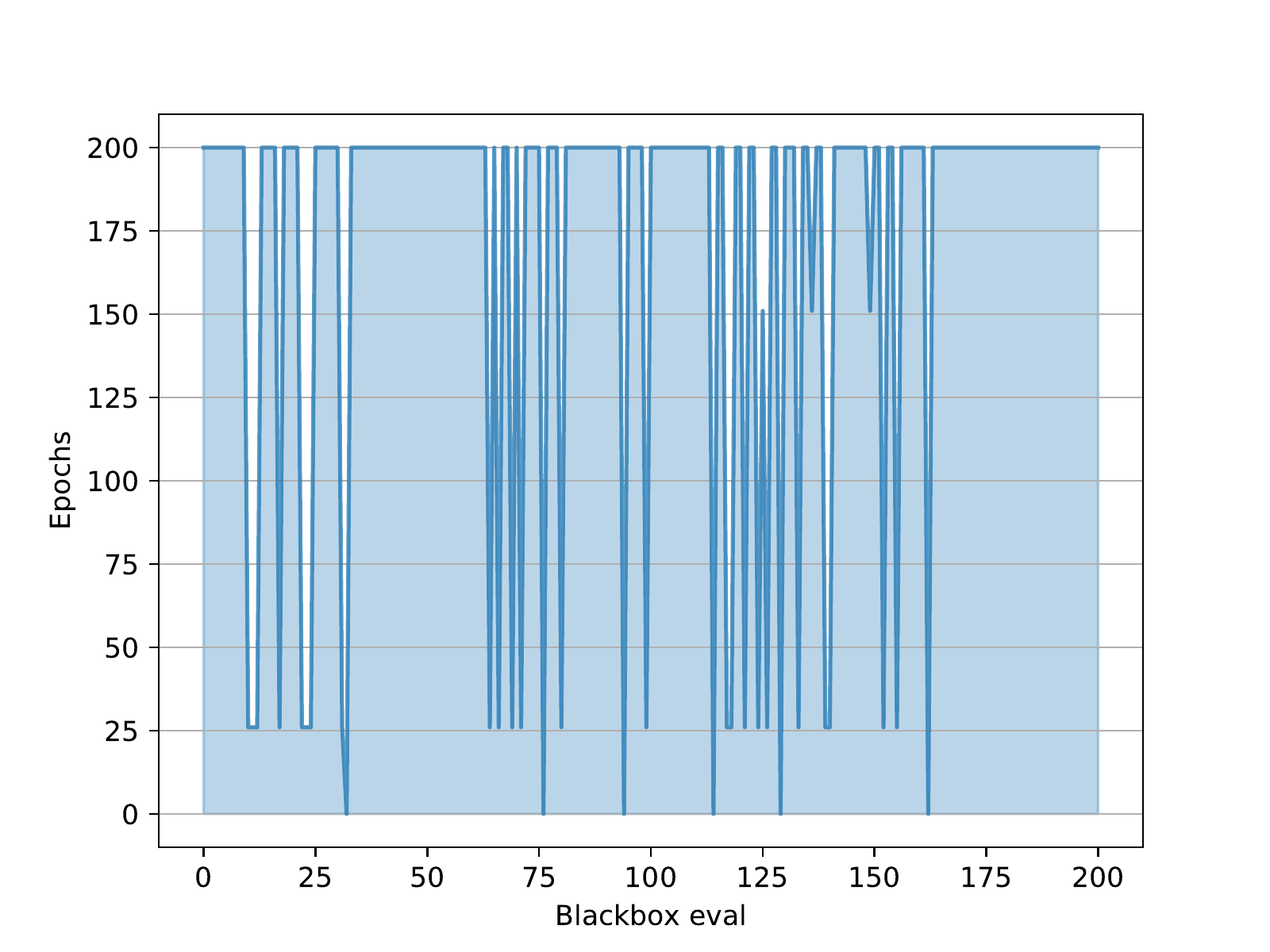}
        \caption{\hypernomad: Stop if the standard deviation of the validation loss is lower than $1e^{-3}$.}
        \label{fig:hn_def_mnist}
    \end{subfigure}
    \hfill
    % \hspace{1mm}
    \begin{subfigure}[t]{0.49\textwidth}
        \centering
        \includegraphics[width=\textwidth]{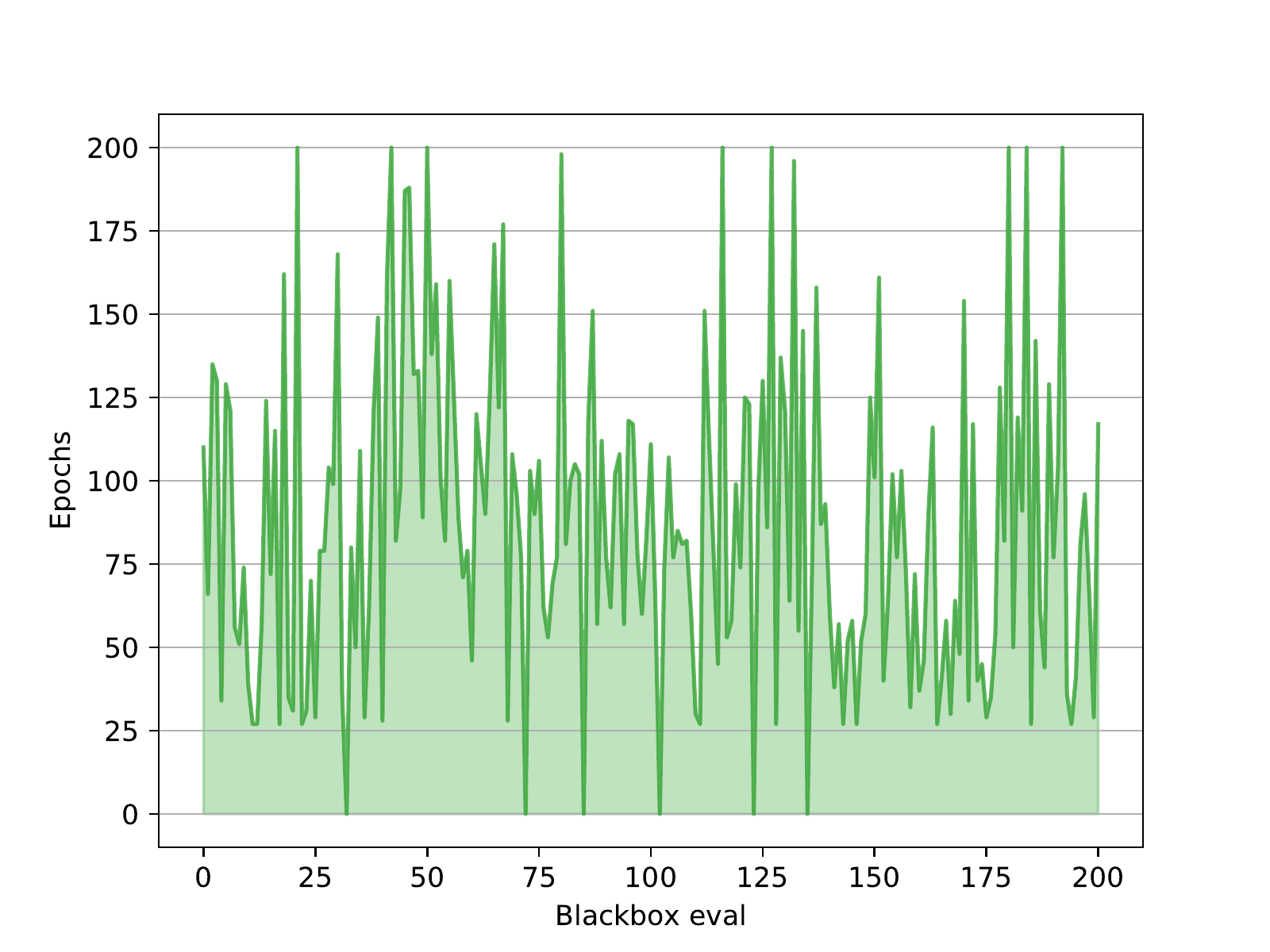}
        \caption{\hypernomad + last success: Stop if the last improvement on the validation score is recorded more than 25 epochs ago.}
        \label{fig:val_acc_fix_mnist}
     \end{subfigure}
     \hfill
    % \hspace{1mm}
    %  \\
    \begin{subfigure}[t]{0.49\textwidth}
        \centering
        \includegraphics[width=\textwidth]{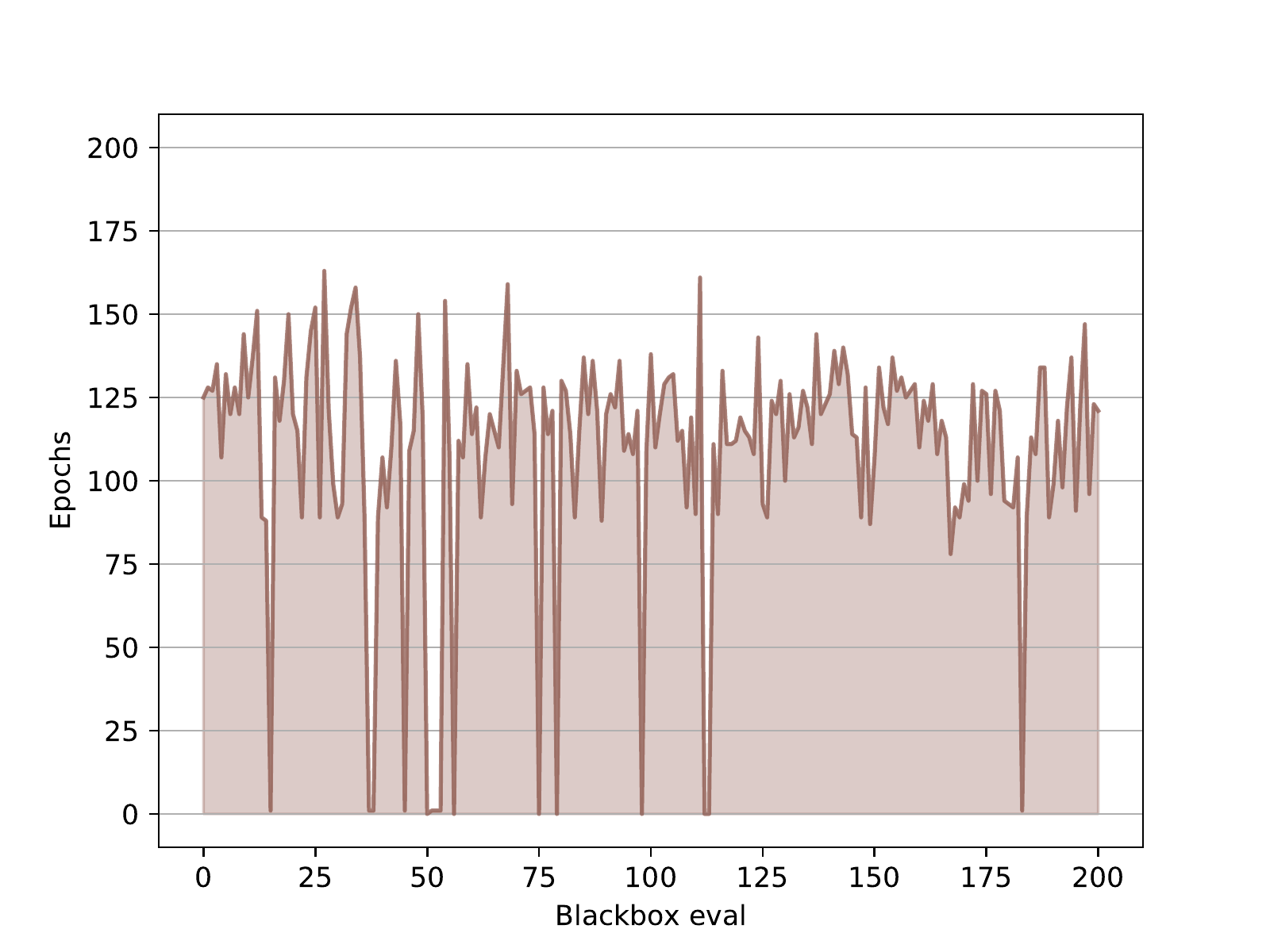}
        \caption{\hypernomad + scheduler: Stop if a plateau is detected with the scheduler.}
        \label{fig:scheduler_def_mnist}
    \end{subfigure}
    \hfill
    % \hspace{1mm}
    \begin{subfigure}[t]{0.49\textwidth}
         \centering
         \includegraphics[width=\textwidth]{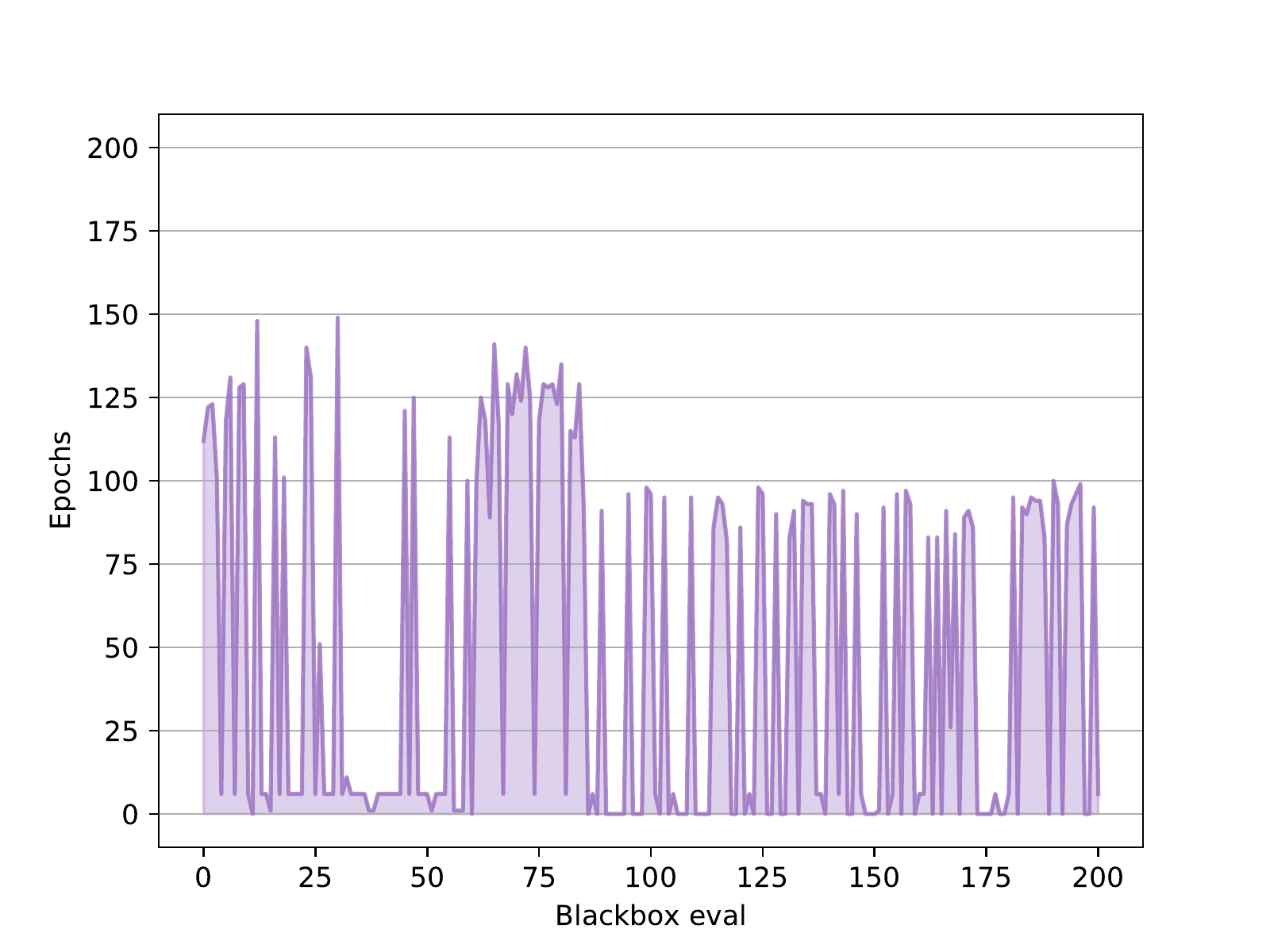}
         \caption{\hypernomad + baseline + scheduler: Stop if a plateau is detected with the scheduler, or if the network scores poorly compared to the baseline. }
         \label{fig:val_acc_var_mnist}
    \end{subfigure}

    \caption{Comparison between the resource consumption of \hypernomad and one variant that adds an early stopping strategy, in terms of number of epochs used to train each visited candidate during a single HPO. The HPO starts from the default point $p_1$ on the MNIST dataset.}
    \label{fig:mnist_def}
\end{figure}

Figures~\ref{fig:mnist_def} displays the number of epochs used in training each of the 200 configurations evaluated during each execution. The original version is the most resource-consuming, with a mean number of epochs per training at 170. All three added early stopping strategies manage to reduce this measure. The mean number of epochs per training drops to 101.4 with the scheduler, 72.5 with the last success criteria, and combining the scheduler with the baseline comparison brings it down to 45.9. Tables~\ref{tab:mnist_def} and~\ref{tab:mnist_pt2} corroborate this conclusion as the proposed approach is the quickest by a factor of 4 in the second example, plus the score of the best solution does not appear to significantly deteriorate compared to the original version of \hypernomad. 
% Test the impact of these stops.
% Analyse the results

\begin{table}[ht]
    \caption{Comparison between four early stopping strategies implemented in \hypernomad in terms of top-1 validation accuracy, overall wall clock time and total number of epochs. The variants are launched on MNIST from two starting configurations, are allowed a total of 200 configuration trials and each network can train during a maximum of 200 epochs.}
    \begin{scriptsize}
    \centering
    \begin{subtable}[b]{0.4\textwidth}
        \caption{Scores of \hypernomad on MNIST with each stopping criteria starting from the default point~$p_1$.}
        \centering
        \hspace*{-8mm}
        \begin{tabular}
        {>{\raggedright}p{0.36\textwidth}>{\centering}p{0.2\textwidth}>{\centering}p{0.15\textwidth}>{\raggedleft\arraybackslash}p{0.15\textwidth}}
        \toprule
        {Early stopping strategy} & Top-1 val. acc. & Wall-clock time (s) & Total epochs \\ \midrule
        % \hline
        \hypernomad & $99.38 \%$  & $305856$ &  $35503$ \\
        Last success & $99.40\%$ & $129600$ & $17534$ \\
        Scheduler  & $99.29\%$ & $170208$ & $21987$\\
        % Validation accuracy & $99.47 \%$ & $0.66$ &  $23867$\\
        Scheduler and~baseline & $\mathbf{99.41\%}$ & $\mathbf{126 144}$ & $\mathbf{9681}$ \\
        \bottomrule
        \end{tabular}
        \label{tab:mnist_def}
    \end{subtable}
    \hspace{5mm} 
    \begin{subtable}[b]{0.4\textwidth}
        \caption{Scores of \hypernomad on MNIST with each stopping criteria starting from the configuration~$p_2$.}
        \centering
        \begin{tabular}{>{\raggedright}p{0.36\textwidth}>{\centering}p{0.2\textwidth}>{\centering}p{0.15\textwidth}>{\raggedleft\arraybackslash}p{0.15\textwidth}}  \toprule
        Early stopping strategy & Top-1 val. acc. & Wall-clock time (s) & Total epochs \\ \midrule
        % \hline
        \hypernomad & $99.43 \%$  & $280800$ & $32712$  \\
        Last success & $\mathbf{99.51\%}$  & $280800$ & $ 11480$ \\
        Scheduler  & $99.29 \%$ & $170208$ & $18580$\\
        Scheduler and~baseline & $99.44\%$ & $\mathbf{64800}$ & $\mathbf{8655}$ \\
        \bottomrule
    \end{tabular}
    \label{tab:mnist_pt2}
    \end{subtable}
    \end{scriptsize}
\end{table}

%---------------------------------------------------------------%
\subsection{Ranking}
%---------------------------------------------------------------%

In addition to early stopping, ranking a new pool of candidates based on the expected performance can accelerate the overall HPO process. In \hypernomad, the evaluations are opportunistic, meaning that the poll step at iteration $k$ stops as soon as a point $p_k \in P_k$ scores better than the incumbent $x_k$. Therefore, the opportunistic strategy prioritizes fast improvements and, when coupled with an adequate ranking tool, has the potential to explore the search space more efficiently and find the best configurations that much quicker.

Static surrogates, or low fidelity estimates, are approximating functions that do not change during the HPO execution, contrary to dynamic surrogates, which are updated to account for the collected information with each iteration. Creating a low fidelity surrogate to estimate the accuracy of a DNN can mean training on a fraction of the full epoch budget or on a subset of the dataset in this particular setting. The ideal static surrogate combines a cheap evaluation with a reliable estimation to correctly rank a pool of candidates. Note that in the context of \hypernomad and \mads, a good static surrogate does not need to produce accurate approximations of the true objective but rather preserves the ranking order between the candidates. 

The effects of ranking the poll points are observed on the MNIST dataset by comparing four static surrogates added to \hypernomad. Two surrogates train the candidates on a low training budget, and two others train on a subset of the dataset. Recall that the proposed early stopping strategy results in a mean epoch budget per blackbox evaluation at around 45. Therefore, a static surrogate must use a lower training budget to qualify as a cheap approximation of the true objective evaluation. Similarly, only a small portion of the dataset can be used for the second set of surrogates. The tested surrogates are summarized in Table~\ref{tab:ranking}.

\begin{figure}[ht]
    \begin{subfigure}[t]{0.495\textwidth}
        \centering
        \includegraphics[width=\textwidth]{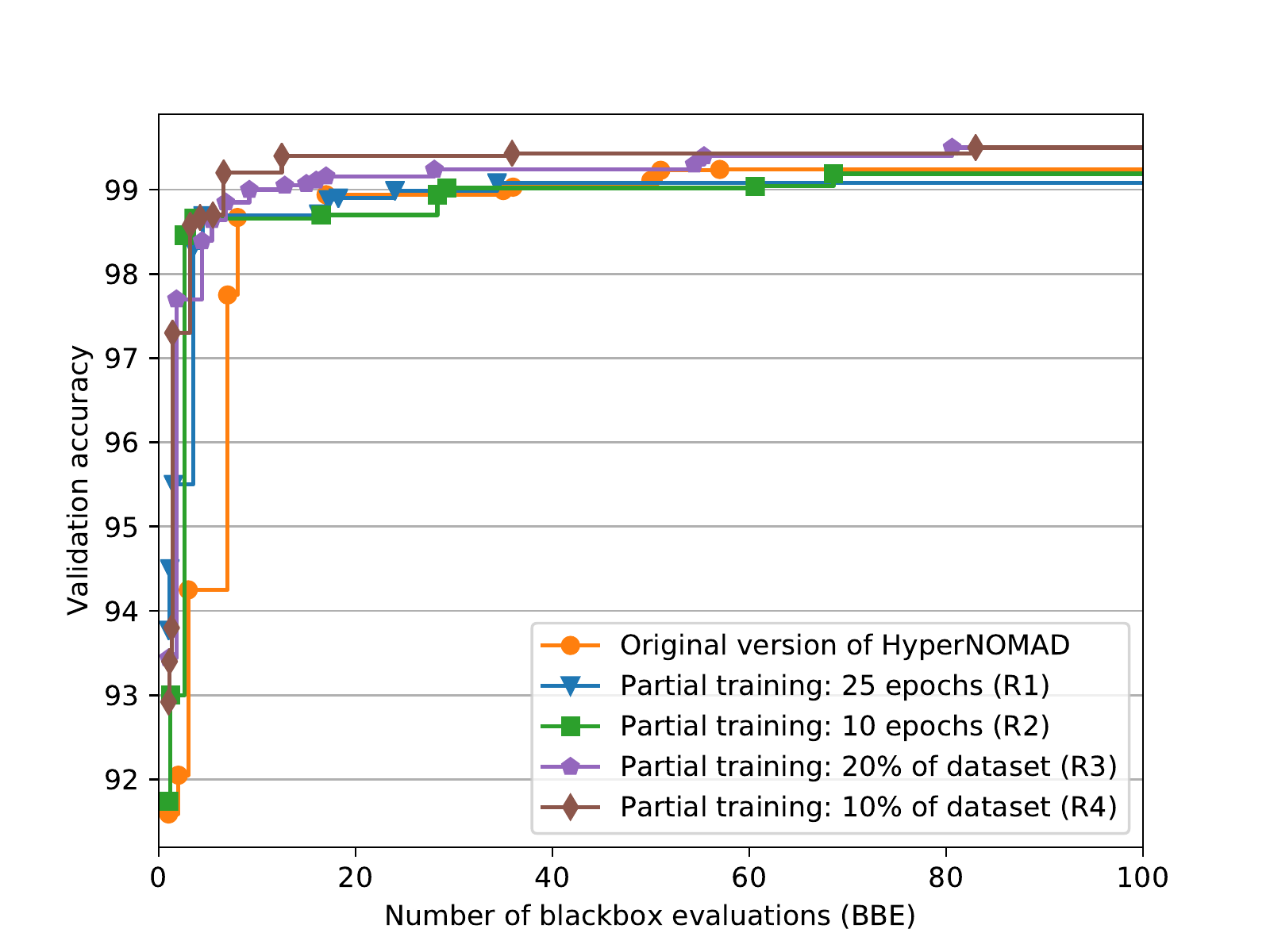}
        \caption{Convergence of each variant in terms of validation accuracy per number of BBE.}
        \label{fig:acc_bb_pt1}
    \end{subfigure}
    \hfill
    \begin{subfigure}[t]{0.495\textwidth}
        \centering
        \includegraphics[width=\textwidth]{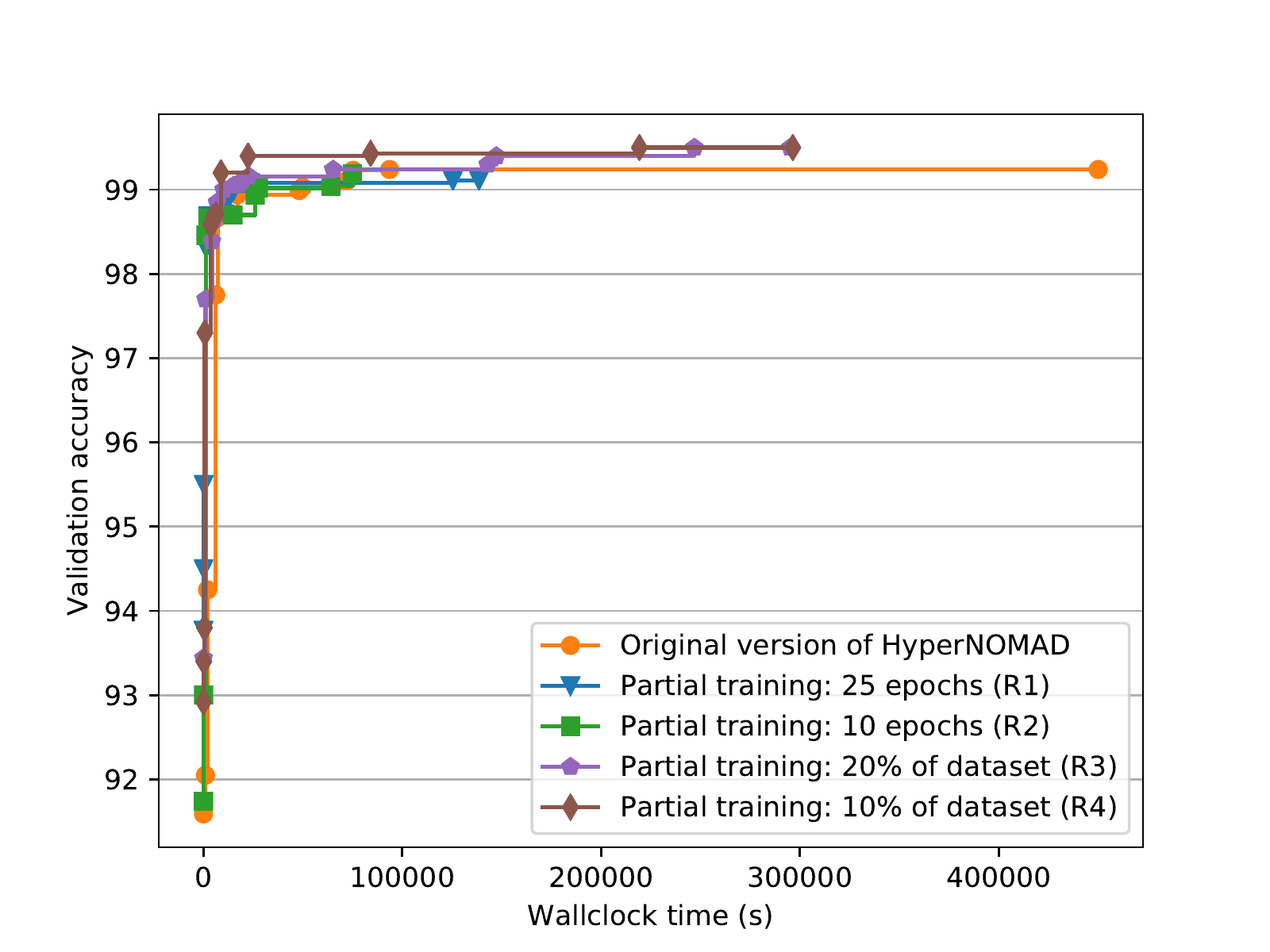}
        \caption{Convergence of each variant in terms of validation accuracy per overall execution time.}
        \label{fig:acc_time_pt1}
     \end{subfigure} 
     
     \caption{Comparison between the original \hypernomad~\cite{hypernomad,hypernomad_paper} and four variants that implement a strategy to sort the new candidates in order to evaluate the most promising first. Strategies $R_1$ and $R_2$ train on a small epoch budget while strategies $R_3$ and $R_4$ train on a subset of the data. The optimization is launched on the MNIST dataset from the default point $p_1$.}
     \label{fig:sorting_sgte}
\end{figure}

\begin{figure}[ht]
    \begin{subfigure}[t]{0.49\textwidth}
        \centering
        \includegraphics[width=\textwidth]{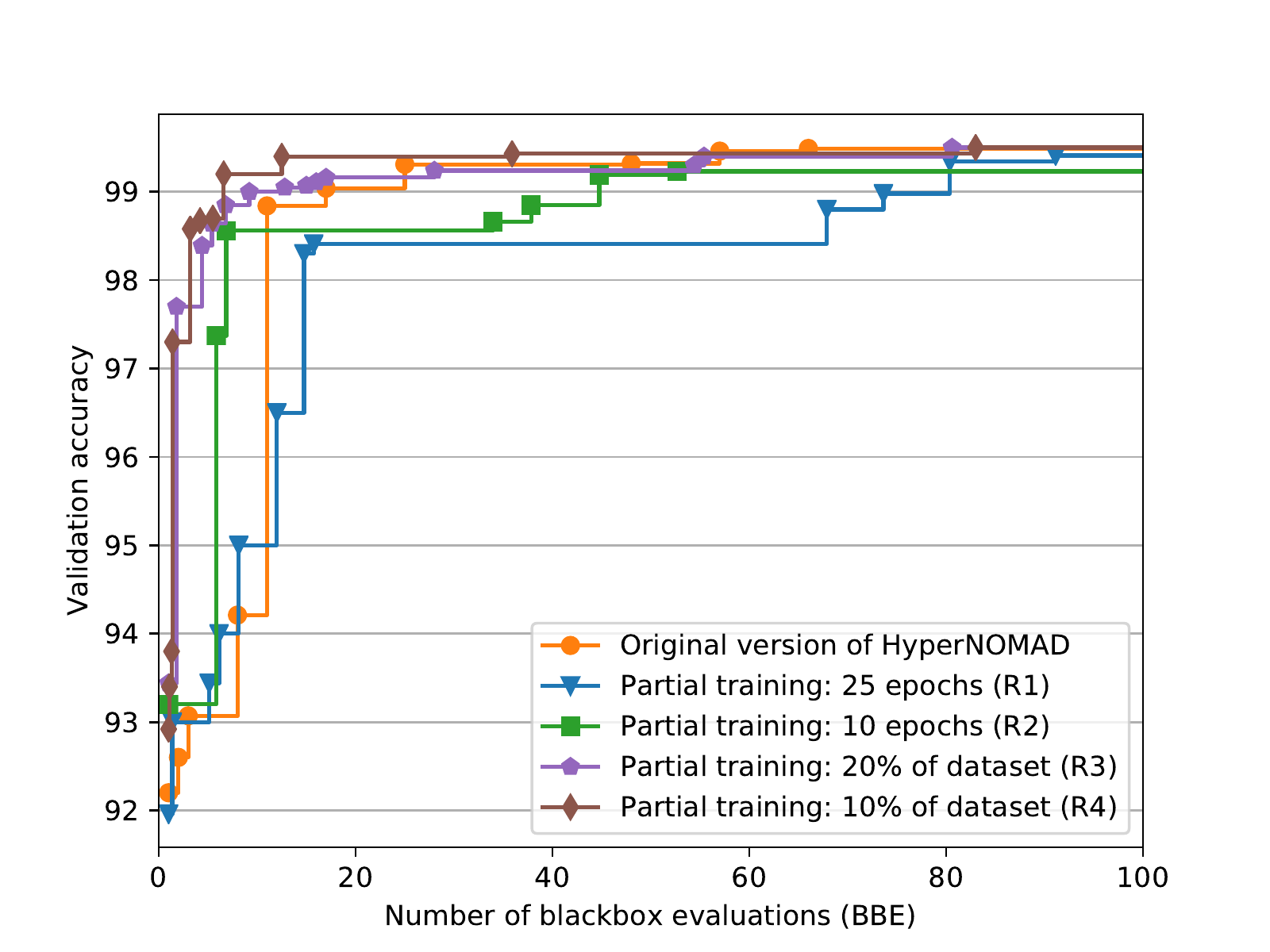}
        \caption{Convergence of each variant in terms of validation accuracy per number of BBE.}
        \label{fig:acc_bb_pt2}
    \end{subfigure}
    \hfill
    \begin{subfigure}[t]{0.49\textwidth}
        \centering
        \includegraphics[width=\textwidth]{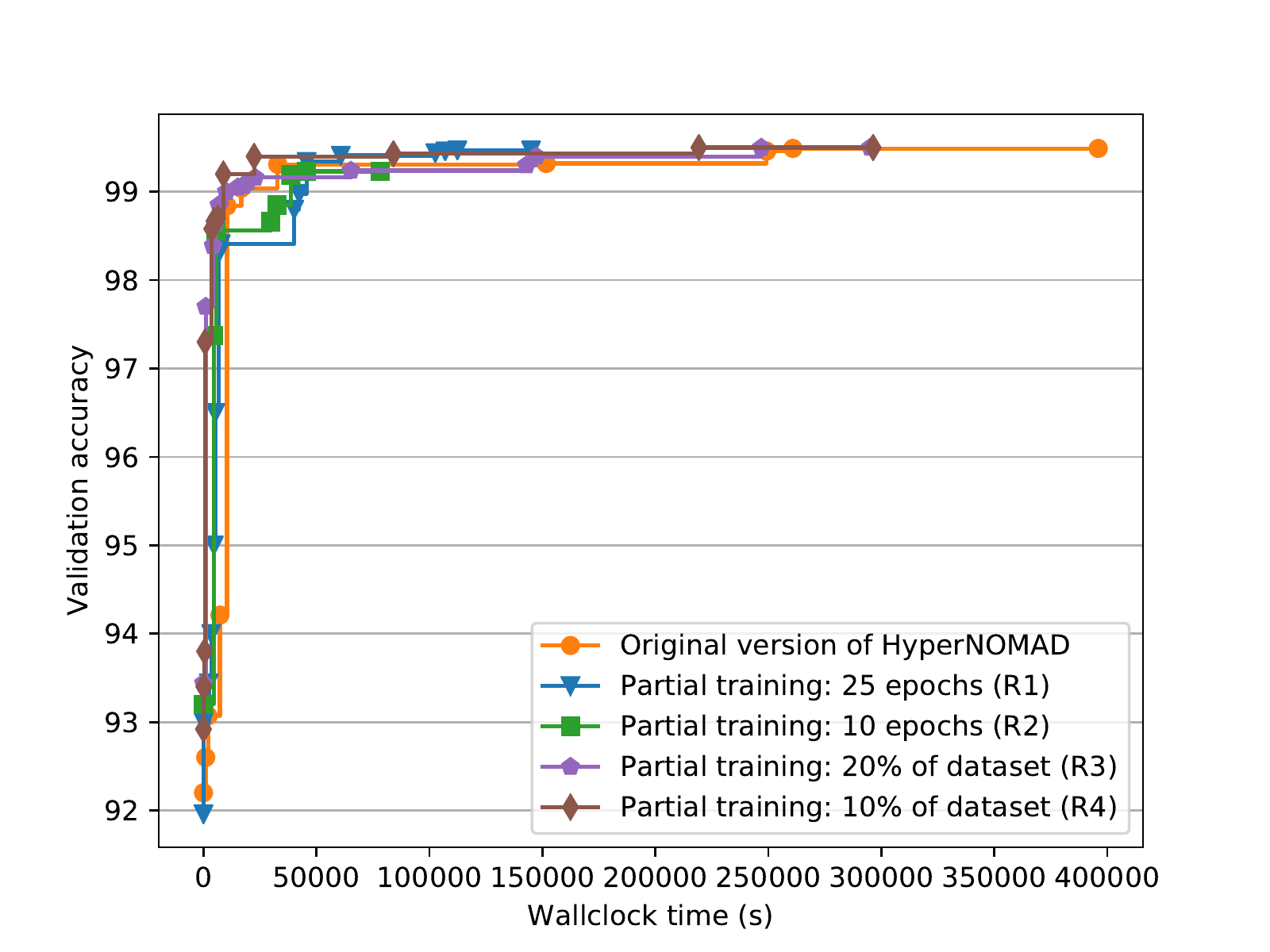}
        \caption{Convergence of each variant in terms of validation accuracy per overall execution time.}
        \label{fig:acc_time_pt2}
     \end{subfigure} 
     
     \caption{Comparison between the original \hypernomad~\cite{hypernomad,hypernomad_paper} and four variants that implement a strategy to sort the new candidates in order to evaluate the most promising first. Strategies $R_1$ and $R_2$ train on a small epoch budget while strategies $R_3$ and $R_4$ train on a subset of the data. The optimization is launched on the MNIST dataset from the default point $p_2$.}
     \label{fig:sorting_sgte2}
\end{figure}

All variants start from the same configurations $p_1$ and $p_2$, are allowed 100 blackbox evaluations with a budget of 200 epochs each, and no early stopping criteria are considered at this stage. Evaluating a surrogate must also account in the number of blackbox evaluations as they add to the resource consumption of such HPO optimization. The cost of each surrogate in terms of blackbox evaluation is provided in Table~\ref{tab:ranking}. Figures~\ref{fig:sorting_sgte} and~\ref{fig:sorting_sgte2} summarize the results regarding the best validation accuracy per number of blackbox evaluations and overall clock-time. Once again, all executions are run on a single Nvidia P100 GPU.

%---------------------------------------------------------------%
\begin{table}[ht]
    \caption{List of the static surrogates tested for ranking the pool candidates. Their evaluation cost in terms of epochs and portion of datasets is compared to a full blackbox evaluation (BBE).}
    \centering
        \begin{tabular}{>{\raggedright}p{0.2\textwidth} >{\centering}p{0.2\textwidth} >{\centering}p{0.2\textwidth}>{\raggedleft\arraybackslash}p{0.2\textwidth}}  \toprule
        % {l|c|c|c }
        Function & Training budget (epochs) & Portion of dataset & Cost ratio to full BBE \\ \midrule
        % \hline
        Objective function & $200$  & $100\%$ &  $100\%$ \\
        Surrogate $R_1$ & $25$ & $100\%$ & $12.5\%$ \\
        Surrogate $R_2$  & $10$ & $100\%$ & $5\%$\\
        Surrogate $R_3$ & $200$ & $20\%$ & $20\%$ \\
        Surrogate $R_4$  & $200$ & $10\%$ & $10\%$\\
        \bottomrule
        \end{tabular}
    \label{tab:ranking}
\end{table}
%---------------------------------------------------------------%

Figures~\ref{fig:acc_bb_pt1} and~\ref{fig:acc_bb_pt2} show the effect of adding ranking surrogate on the overall convergence of \hypernomad. In Figure~\ref{fig:acc_bb_pt1}, all variants are relatively equivalent and surpass \hypernomad in the early stages of the execution. This tendency shifts however after around 10 blackbox evaluations when $R_3$ and $R_4$ score better quality solutions faster than the other variants. In fact, $R_1$ and $R_2$ have limited benefits to \hypernomad since the latter appears to have an equivalent performance in this case. This observation is more obvious in Figure~\ref{fig:acc_bb_pt2} where adding $R_1$ or $R_2$ significantly slow down the execution of \hypernomad.

Figures~\ref{fig:acc_time_pt1} and~\ref{fig:acc_time_pt2} compare the convergence of each algorithm in term of overall wall-clock time. The benefit of adding a sorting surrogate is apparent in both executions as \hypernomad is the variant that takes the longest to evaluate the 100 blackbox trials. This is a coherent observation with the fact that every one of the 100 evaluations is equivalent to fully training a network on the entire dataset. Also, in both series of tests, surrogates $R_1$ and $R_2$ are the fastest and often terminate the 100 blackbox evaluation budget before all the other variants. It appears that training on such low epoch budgets, even on the full dataset, is faster than training on $10\%$ of the dataset during 200 epochs. 
However, this observation is expected to change when early stopping is integrated again, which stops $R_3$ and $R_4$ from using the full 200 epochs in their evaluation. 
Plus, since $R_1$ and $R_2$ harm the best configuration score, both these strategies are discarded, and it is $R_4$ that is retained.

%---------------------------------------------------------------%
\section{Testing}\label{sec:results}
%---------------------------------------------------------------%

This section incorporates both aspects of Section~\ref{sec:surrogates} into the \hypernomad framework. Early stopping is triggered from a scheduler or by comparison with the baseline network; and ranking a new pool of candidates is achieved through training each candidate on 10\% of the dataset. 

The tests are now launched on the CIFAR-10 dataset, are allowed 100 blackbox evaluations with a maximum training budget of 200 epochs per configuration. This time, each execution is launched on two Nvidia P100 GPUs. One series of tests starts from the default configuration of \hypernomad noted $p_1$, and the second starts from a network with a 5 convolutional layers and 1 fully connected layer, noted $p_3$, which is equivalent to 36 hyperparameters. 

\begin{figure}[ht]
    \begin{subfigure}[b]{0.49\textwidth}
        \centering
        \includegraphics[width=\textwidth]{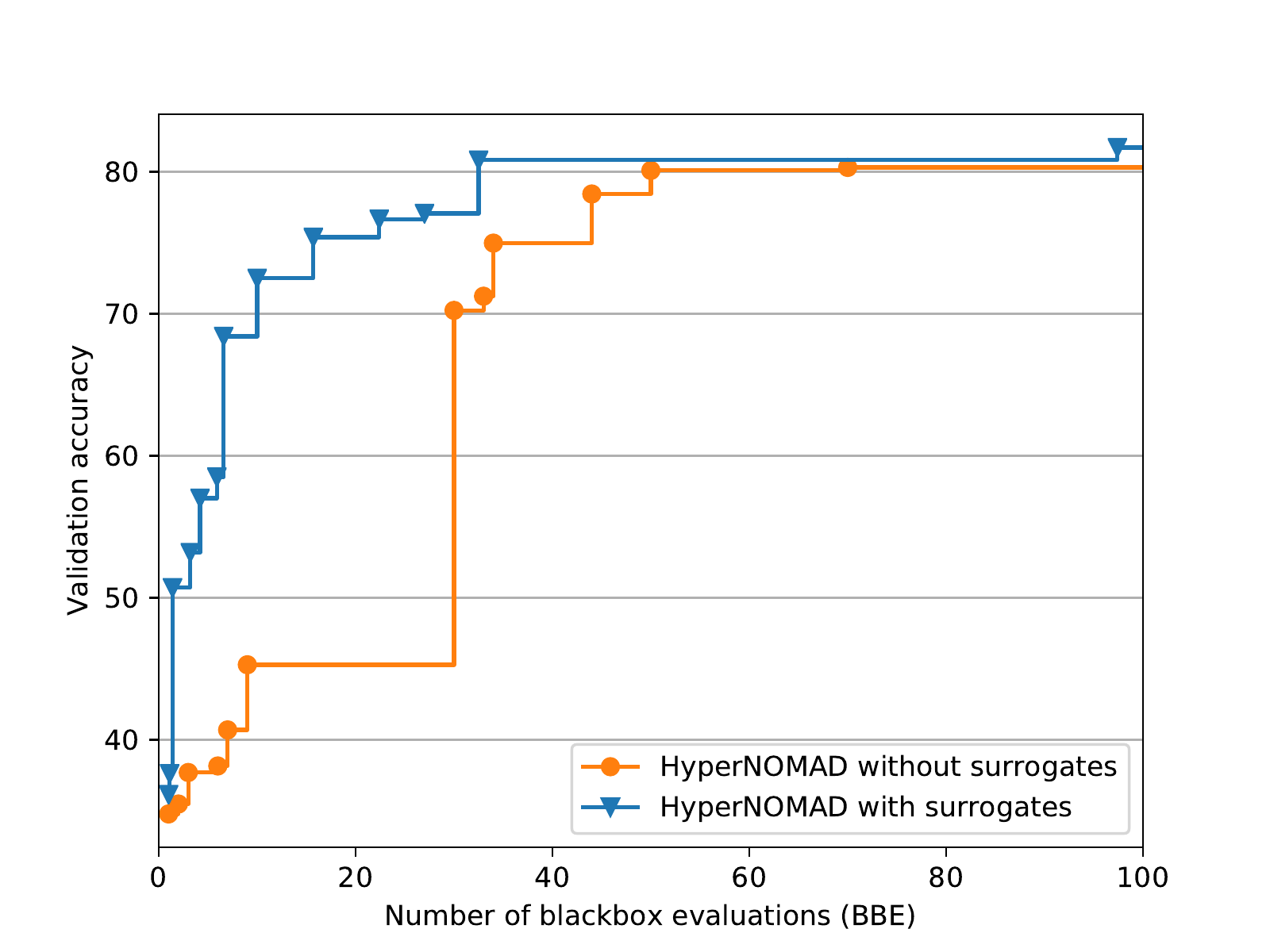}
        \caption{Convergence in terms of validation accuracy per number of BBE.}
        \label{fig:acc_bb_pt1_cifar}
    \end{subfigure}
    \hfill
    \begin{subfigure}[b]{0.49\textwidth}
        \centering
        \includegraphics[width=\textwidth]{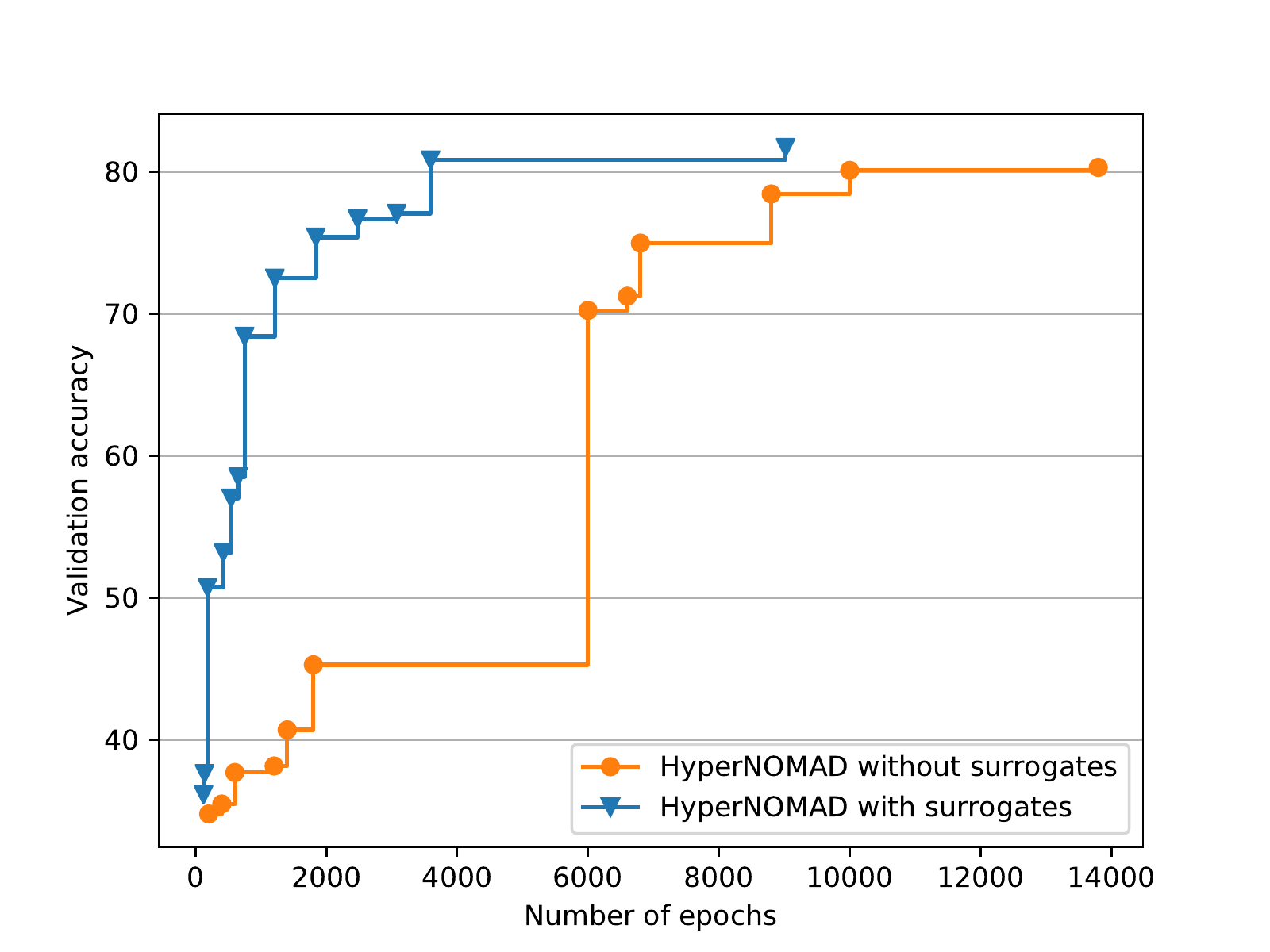}
        \caption{Convergence in terms of validation accuracy per number of epochs.}
        \label{fig:acc_epochs_pt1_cifar}
    \end{subfigure}
    \hfill
    \begin{subfigure}[b]{0.49\textwidth}
        \centering
        \includegraphics[width=\textwidth]{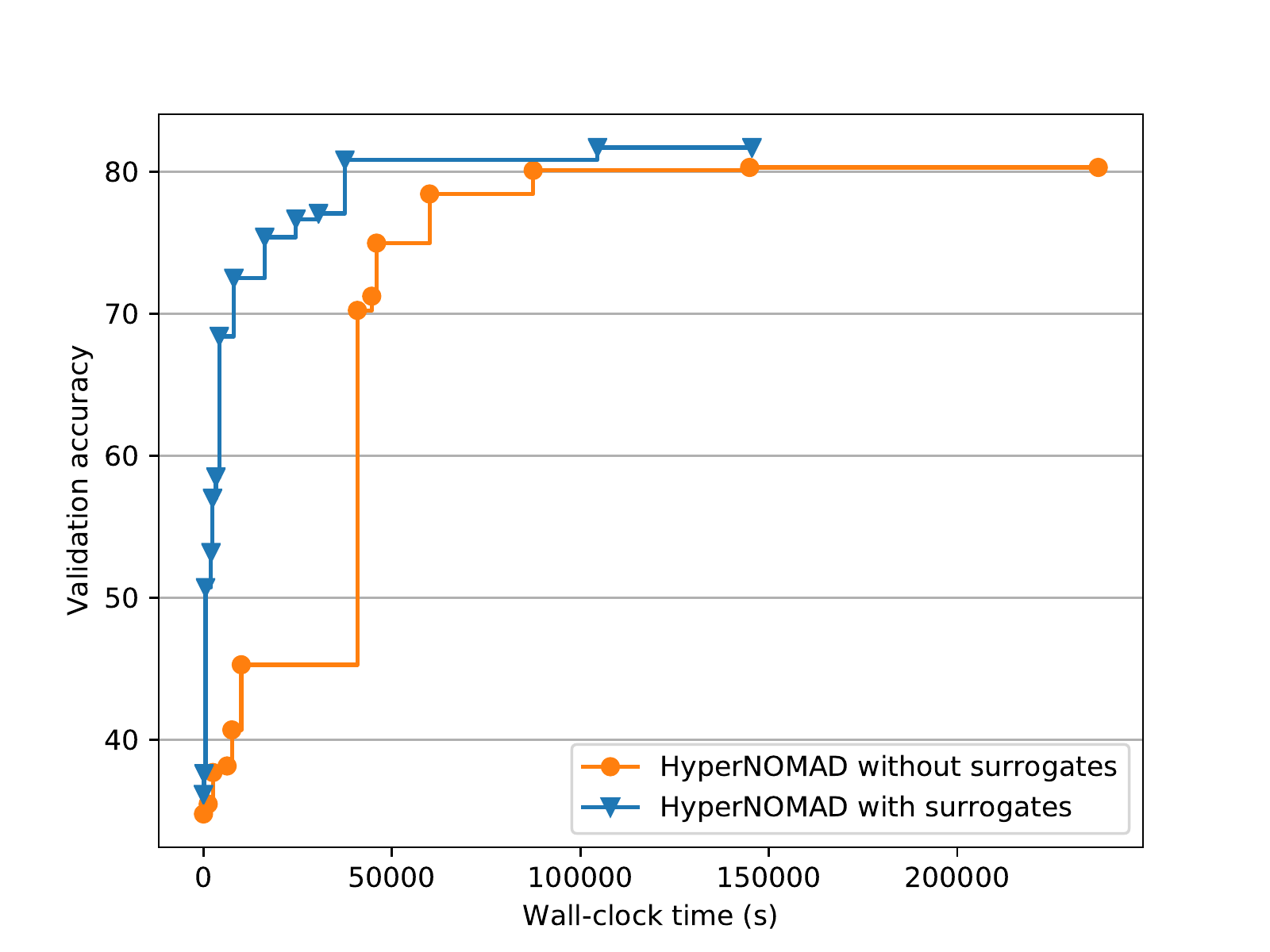}
        \caption{Convergence in terms of validation accuracy per execution time.}
        \label{fig:acc_time_pt1_cifar}
     \end{subfigure} 
     
     \caption{Comparison between \hypernomad with and without surrogates on the CIFAR-10 dataset, starting from the default settings $p_1$.}
     \label{fig:cifar_pt1}
\end{figure}

\begin{figure}[ht]
    \begin{subfigure}[b]{0.49\textwidth}
        \centering
        \includegraphics[width=\textwidth]{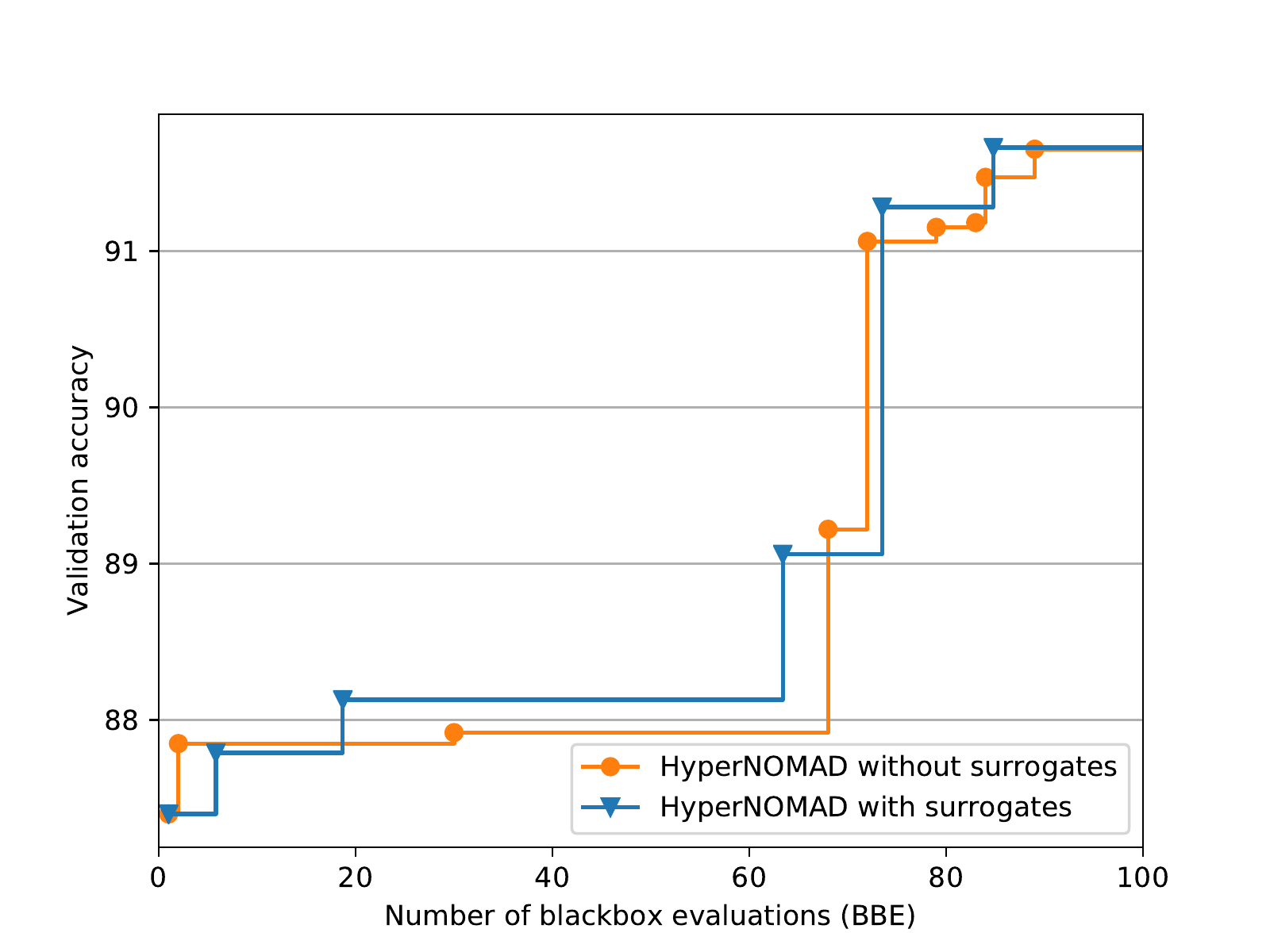}
        \caption{Convergence in terms of validation accuracy per number of BBE.}
        \label{fig:acc_bb_pt3_cifar}
    \end{subfigure}
    \hfill
    \begin{subfigure}[b]{0.49\textwidth}
        \centering
        \includegraphics[width=\textwidth]{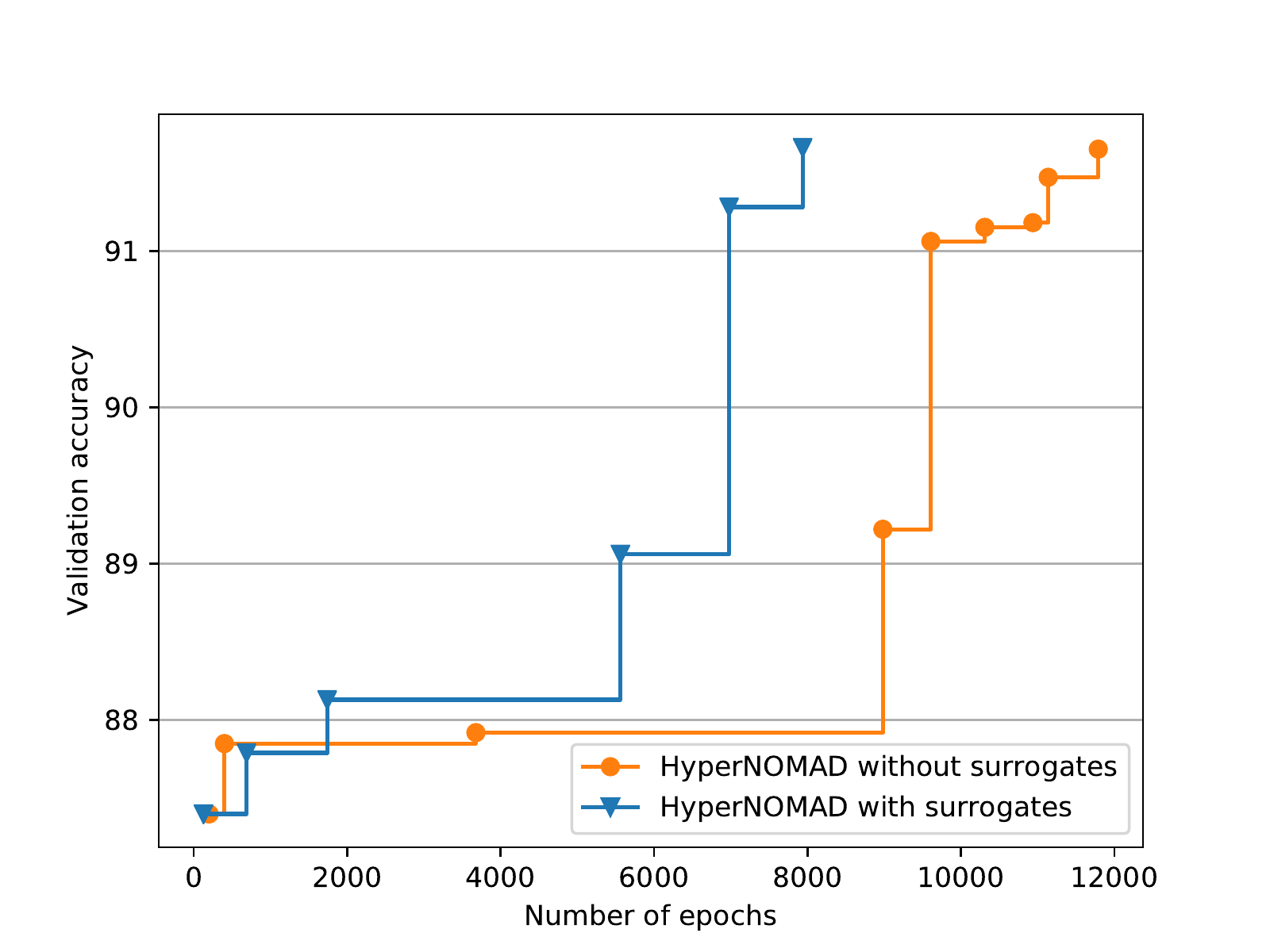}
        \caption{Convergence in terms of validation accuracy per number of epochs.}
        \label{fig:acc_epochs_pt3_cifar}
    \end{subfigure}
    \hfill
    \begin{subfigure}[b]{0.49\textwidth}
        \centering
        \includegraphics[width=\textwidth]{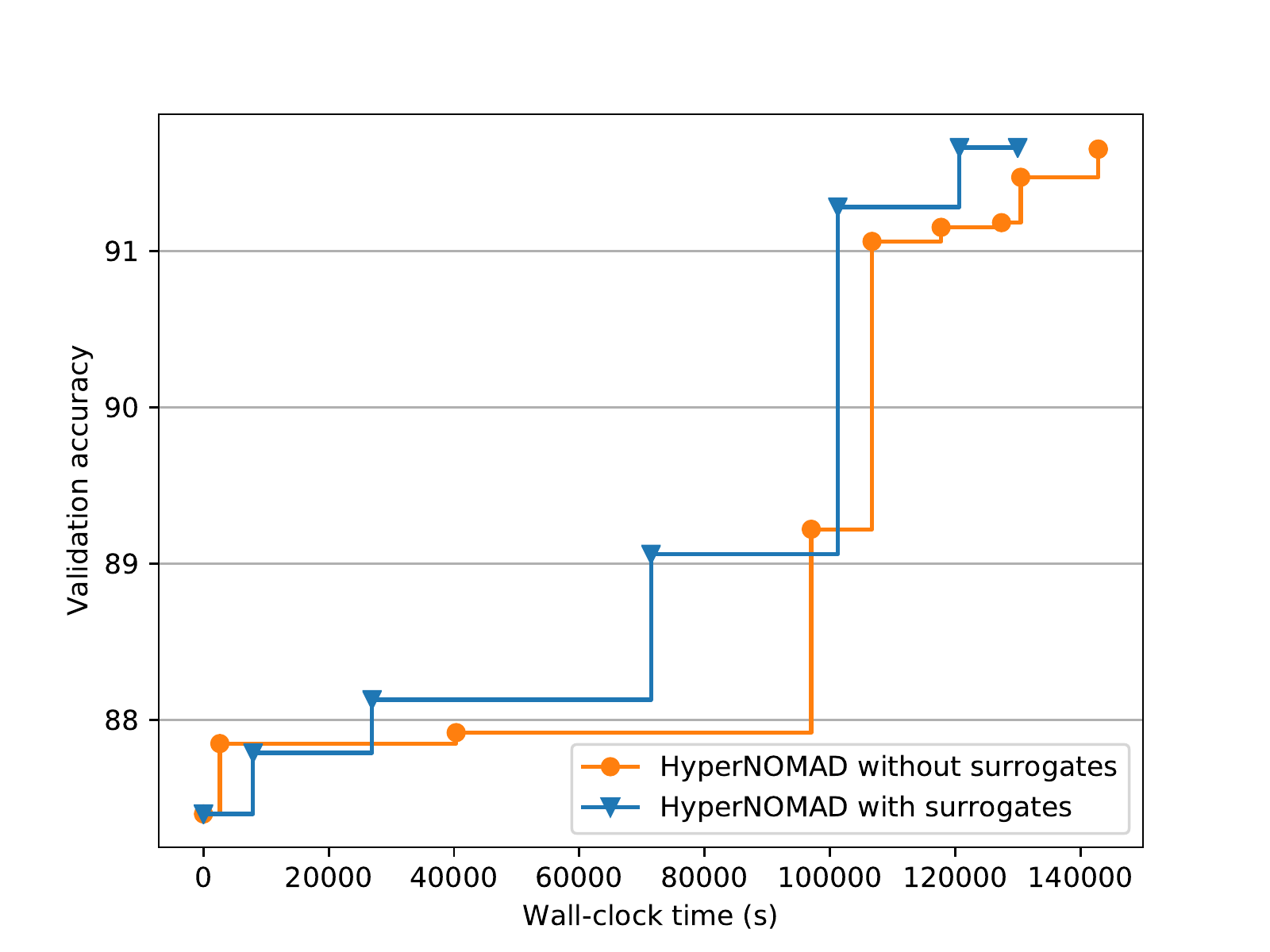}
        \caption{Convergence in terms of validation accuracy per execution time.}
        \label{fig:acc_time_pt3_cifar}
     \end{subfigure} 
     
     \caption{Comparison between \hypernomad with and without surrogates on the CIFAR-10 dataset, starting from the default settings $p_3$.}
     \label{fig:cifar_pt3}
\end{figure}

Figure~\ref{fig:cifar_pt1} summarizes the results from the comparison between the original version of \hypernomad and the one with the proposed enhancements when launched from the default settings $p_1$. The benefits of early stopping and ranking the new candidates are apparent, with the new version scoring better than \hypernomad in terms of best validation score per blackbox evaluation budget as seen in Figure~\ref{fig:acc_bb_pt1_cifar}. Figure~\ref{fig:acc_epochs_pt1_cifar} also shows that adding surrogates saves training resources in terms of total number of epochs to get to better quality solution, and Figure~\ref{fig:acc_time_pt1_cifar} corroborates this tendency when comparing the overall wall-clock time to deplete all 100 blackbox evaluations.
The second test starts from configuration $p_3$. Figure~\ref{fig:acc_bb_pt3_cifar} shows that the version with the surrogates is slightly ahead when comparing the best validation score per blackbox evaluation budget. Figure~\ref{fig:acc_epochs_pt3_cifar} illustrates a gain of 27\% in total epoch budget, yet this gain is not translated in terms of overall wall-clock time as shown in Figure~\ref{fig:acc_time_pt3_cifar}. The difference is believed to be caused by the communication between the two GPUs used at once.

%------------------------------------------------%
\section{Discussion}
%------------------------------------------------%
This work proposes a solution to speed up the time and resource-consuming process of optimizing the hyperparameters of deep neural networks, based on incorporating two strategies based on static surrogates into the \hypernomad framework. The first defines a new early stopping strategy that quickly and effectively interrupts poorly performing networks. The second allows ranking a pool of new candidates to detect and then evaluate the most promising first. Both techniques are shown individually and collectively to improve on the quality of solution and on the amount of computational resources needed. 

%------------------------------------------------%
\section*{Acknowledgments}
%------------------------------------------------%

This work is in part supported by the NSERC Alliance grant 544900-19 in collaboration with Huawei-Canada.

%------------------------------------------------%
% References
%------------------------------------------------%
%\bibliographystyle{plain}
%\bibliography{bibliography}
%------------------------------------------------%

%---------------------------------------------------------------%
\end{document}